# A global DC branch model incorporating power system flexibility

F. M. Gatta, A. Geri, S. Lauria, M. Maccioni, L. Nati

Dept. of Astronautics, Electrical and Energy Engineering (DIAEE), “Sapienza” University of Rome, Rome, Italy

**Abstract**

In this letter we propose a generalized branch model to be used in DC optimal power flow (DCOPF) applications. Besides AC lines and transformers, the formulation allows for representing variable susceptance branches, phase shifting transformers, HVDC lines, zero impedance lines and open branches. The possibility to model branches with concurrently variable susceptance and controllable phase shift angles is also provided. The model is suited for use in DCOPF formulations aimed at the optimization of remedial actions so as to exploit power system flexibility; applications to small-, medium- and large-scale systems are presented to this purpose.

**Keywords:** DC optimal power flow, FACTS devices, HVDC lines, phase shifting transformers, transmission switching, zero impedance lines.

**Nomenclature**

Quantities in bold characters refer to vectors and matrices.

*Indices and Sets*

$k$, $n$, $t$: indices of branches (running from 1 to $N_K$), of nodes (running from 1 to $N_N$) and of time periods (running from 1 to $N_T$) respectively.

*Parameters*

$\boldsymbol{f_{max}}$: $N_T \times N_K$ matrix of maximum branch capacities.

$\boldsymbol{p_c}$: The same of $\boldsymbol{f_{max}}$, but capacities of branches different than HVDC lines and zero impedance lines are set to 0.

$\boldsymbol{\varphi_{min}}$, $\boldsymbol{\varphi_{max}}$: $N_T \times N_K$ matrices of minimum and maximum shift angles.

$\boldsymbol{\Delta\delta_{max}}$: $N_T \times N_K$ matrix of maximum voltage angle differences acceptable at the branch terminals.

$\boldsymbol{B_{min}}$, $\boldsymbol{B_{max}}$: $N_T \times N_K$ matrices of minimum and maximum branch susceptances.

***Cft***: $N_K \times N_N$ connection matrix.

*M*: Large enough number.

*Real Variables*

***f***: $N_T \times N_K$ matrix of overall power flows at branch terminals.

**φ**: $N_T \times N_K$ matrix of shift angles.

$\boldsymbol{p_{br}}$: $N_T \times N_K$ matrix of power flows through branches.

$\boldsymbol{p_{top}}$: $N_T \times N_K$ matrix of power injections at branch terminals modeling open branches.

***nex***: $N_T \times N_N$ matrix of net power injections at nodes.

**δ**: $N_T \times N_N$ matrix of voltage phase angles.

*Integer Variables*

***tra***: $N_T \times N_K$ matrix related to topological changes. $tra(t,k)$ is 1 if branch *k* is open at time *t*, 0 otherwise.

***z***: $N_T \times N_K$ matrix related to power flow direction through a branch with variable susceptance. $z(t,k)$ is 0 if the flow on branch *k* at time *t* is oriented according to the positive direction, otherwise is 1.

## 1 Introduction

The pursuit of optimal power system operation is normally supported by DC optimal power flow (DCOPF) procedures, for which a linearized DC representation of each network component is required. Problems such as network constrained unit commitment, economic dispatch and congestion management can be addressed by controlling power flows across the grid to reduce operational costs and enhancing the power system utilization via appropriate remedial actions. These can be changes to network topology (e.g. transmission switching, substation reconfiguration, bypass of network devices), variations of the HVDC lines' operating points, variations of flexible AC transmission system (FACTS) setpoints and variations of phase shifting transformers (PSTs) taps. Clearly, all the above-mentioned remedial actions are interdependent and an appropriate modeling framework, which can include all of them, is useful to enhance the exploitation of the flexibility of the power system.

Literature addressing the optimization of a transmission network by means of remedial actions is wide and paper are numerous. DCOPF formulations including transmission switching (TS) have been proposed: [1] and [2] use a $B\theta$ formulation, whereas [3] and [4] employ a reduced shift factor formulation. Substation reconfiguration is performed by switching in/off bus-bar couplers, each represented by a zero impedance line (ZILs): in [5], [6] a model representing a substation with only one ZIL is presented, whereas [7] develops a general substation reconfiguration model in order to simulate all substation bus-bar layouts. The inclusion of PSTs in DCOPF formulations is presented in [8]-[10], whereas the further addition of HVDC lines is reported in [9]. Several representations of series compensation devices (TSCS) are available, for instance [8], which linearizes the nonlinear problem by a two-stage approach with Benders decomposition, or [11 ex 12], which employs a mixed integer linear programming (MILP), resolved in a two-stage LP. Reference [12] expands the FACTS formulation in [11] showing the interdependence of TS and series compensation devices, though with an approximated TSCS model with an a priori determined power flow direction.

While the mentioned literature shows a significant correlation between the power system benefit and remedial actions, none of the previous works consider all the flexibility at disposal, being limited to the analyzed subset of the remedial actions. For example, [1]-[4] model only the possibility to switch on and off transmission lines and, likewise, [5]-[7] limit their analysis to substation reconfiguration through ZIL switching. Furthermore, [9] and [12] point out the interdependence of different kind of remedial actions and propose formulations to address their specific issue ([9] combines the use of HVDC and PST, whereas [12] TS and TSCS).

In this paper we overcome the limitation of the previous works proposing an exact modeling in the DC framework which, besides allowing an easy and fully vectorized implementation, makes use of a unified branch representation able to simultaneously include all remedial actions separately considered in the referenced papers. Moreover, the model allows the simulation of FACTS devices able to control, concurrently or selectively, internal susceptance and phase shift angle (VSSA in the rest of the paper), which to our knowledge has never been proposed in precedent works. Another particular feature of the proposed formulation is the possibility to disconnect any kind of branch, and to bypass any device by paralleling it with a ZIL. The model is envisioned as a building block in any DCOPF formulations (e.g. Economic Dispatch, Unit Commitment, Market Clearing, Transmission Expansion Planning, etc..) making use of the whole flexibility available in the grid.

## 2 Mathematical formulation

In a DC power flow formulation, the power flow $f(t,k)$ at time $t$ through a branch $k$ depends on the branch susceptance value, which in principle may also be variable if TSCS devices are represented, and on the voltage angle difference at branch terminals, possibly adding an imposed shift angle value. Moreover, the branch susceptance is not defined in HVDC lines and ZILs. The power flow on a HVDC line is controlled by the converter stations, whereas in a ZIL (for which the voltage angle difference at its terminals is zero) the power flow is a direct consequence of the power balance at its terminal nodes. Lastly, if the branch is disconnected, $f(t,k)$ must be zero. The generic branch $k$ connecting the buses $n_{from}$ and $n_{to}$ is simulated by two controlled generators, as shown in Fig. 1, where $p_{br}(t,k)$ represents the power flow in the branch and $p_{top}(t,k)$ models the power injection required to simulate a disconnected branch, i.e. $tra(t,k)$=1. The MILP model is provided below.

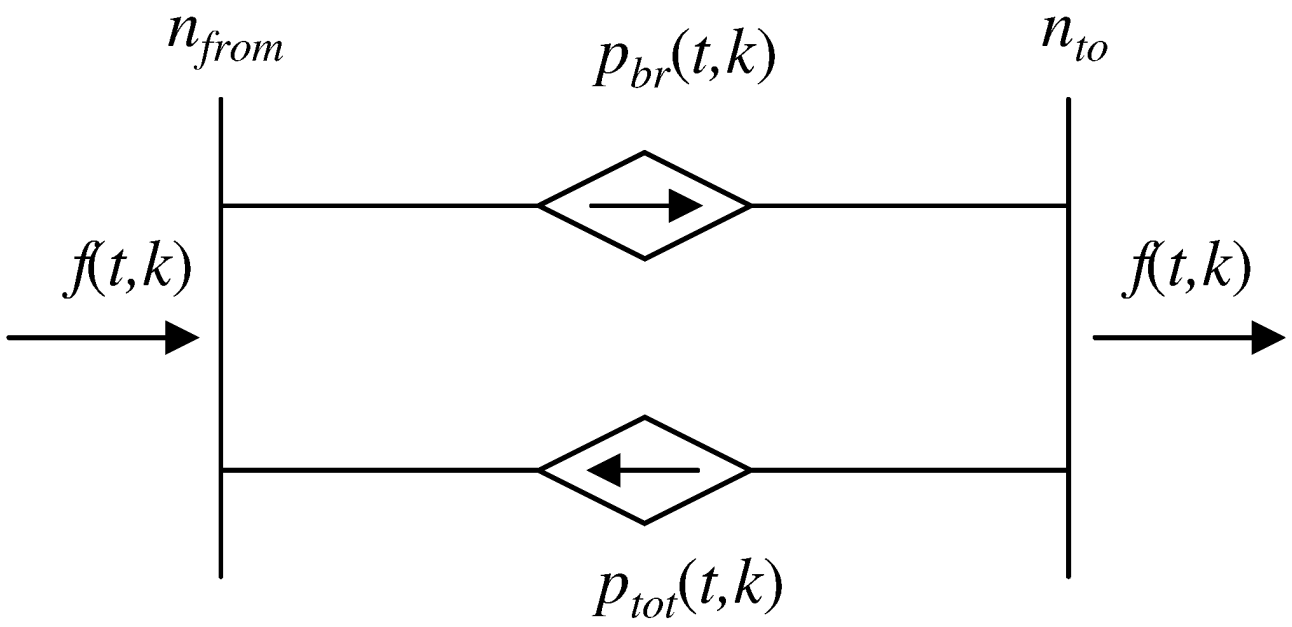

**Fig. 1** – Global branch model in the linearized DC representation.

$$\mathbf{nex} = \mathbf{f} \cdot \mathbf{Cft} \tag{1}$$

$$\mathbf{f} = \mathbf{p_{br}} - \mathbf{p_{top}} \tag{2}$$

$$-M \cdot z - (\mathbf{1} - \mathbf{tra}) \circ \mathbf{p_c} + \mathbf{B_{min}} \circ \left(\delta \cdot \mathbf{Cft}^T + \varphi\right) \leq \mathbf{p_{br}} \leq M \cdot z + (\mathbf{1} - \mathbf{tra}) \circ \mathbf{p_c} + \mathbf{B_{max}} \circ \left(\delta \cdot \mathbf{Cft}^T + \varphi\right) \tag{3}$$

$$-M \cdot (\mathbf{1} - z) - (\mathbf{1} - \mathbf{tra}) \circ \mathbf{p_c} + \mathbf{B_{max}} \circ \left(\delta \cdot \mathbf{Cft}^T + \varphi\right) \leq \mathbf{p_{br}} \leq M \cdot (\mathbf{1} - z) + (\mathbf{1} - \mathbf{tra}) \circ \mathbf{p_c} + \mathbf{B_{min}} \circ \left(\delta \cdot \mathbf{Cft}^T + \varphi\right) \tag{4}$$

$$-M\cdot \boldsymbol{tra}-\Delta\delta_{\boldsymbol{max}}\leq \delta\cdot \boldsymbol{Cft}^{T}\leq \Delta\delta_{\boldsymbol{max}}+M\cdot \boldsymbol{tra} \tag{5}$$

$$(\boldsymbol{1}-\boldsymbol{tra})\circ \varphi_{\boldsymbol{min}}\leq \varphi\leq (\boldsymbol{1}-\boldsymbol{tra})\circ \varphi_{\boldsymbol{max}} \tag{6}$$

$$-(\boldsymbol{1}-\boldsymbol{tra})\circ \boldsymbol{f}_{\boldsymbol{max}}\leq \boldsymbol{f}\leq (\boldsymbol{1}-\boldsymbol{tra})\circ \boldsymbol{f}_{\boldsymbol{max}} \tag{7}$$

$$-M\cdot \boldsymbol{tra}\leq \boldsymbol{p}_{\boldsymbol{top}}\leq M\cdot \boldsymbol{tra} \tag{8}$$

Symbol ◦ refers to Hadamard product and **0** and **1** are $N_T\times N_K$ arrays of zeros and ones, respectively. Equation (1) represent the Kirchhoff's current law (KCL) for the DC power flow representation of the grid. Equations (2)-(8) represent the unified branch model implementation in DC optimal power flow applications. Equation (2) states the power flow through branch terminals, as depicted in Fig. 1. Equations (3) and (4) define the core mathematical conditions of the generalized branch model. For a VSSA branch $k$ with susceptance $B(t,k)$ variable in the range $[B_{min}(t,k), B_{max}(t,k)]$ and with an additional phase shift angle $\varphi(t,k)$ variable in the range $[\varphi_{min}(t,k), \varphi_{max}(t,k)]$, $p_c(t,k)$ is set to zero and binary variable $z(t,k)$ dictates that only one constraint between (3) and (4) is bounding, since the other one collapses to $-M<p_{br}(t,k)<M$. This logical condition is needed to model a branch with variable susceptance with the possibility of both a positive or negative flow through it. Instead, if the branch has a fixed susceptance, i.e. $B_{min}(t,k)=B_{max}(t,k)=B(t,k)$, both (3) and (4) impose $p_{br}(t,k)=B(t,k)\cdot[\delta(t,n_{from})-\delta(t,n_{to})+\varphi(t,k)]$ when bounding, making the $z(t,k)$ value immaterial ($z(t,k)$ may be posed equal either to 0 or 1, indifferently). It is worth noting that VSSA may also represent traditional lines, transformers, PST and TSCS by properly setting lower and upper bounds $[B_{min}(t,k), B_{max}(t,k)]$ and $[\varphi_{min}(t,k), \varphi_{max}(t,k)]$. For HVDC lines and ZILs $B(t,k)$ is not defined, so posing $B_{min}(t,k)=B_{max}(t,k)=0$ (3) and (4) impose $-p_c(t,k)+tra(t,k)\cdot p_c(t,k)\leq p_{br}(t,k)\leq p_c(t,k)-tra(t,k)\cdot p_c(t,k)$ and again $z(t,k)$ may be set indifferently either to 0 or 1. In such a case, $tra(t,k)$ is 0 for a connected branch and so $p_{br}(t,k)$ is limited by the maximum branch capacity $p_c(t,k)$; for a disconnected branch, $tra(t,k)$ is 1 and $p_{br}(t,k)$ is 0. Equation (5) imposes a maximum allowable phase angle difference at branch terminals, which is constrained to zero in case of ZILs; (6) indicates the limits on the additional shift angle which can be imposed at the branch terminals, set to 0 for all the branches which do not encompass this feature. Equation (7) limits $f(t,k)$ through branch terminals between the maximum branch capacity in case $tra(t,k)=0$, otherwise imposes $f(t,k)=0$. Lastly, (8) describes the logical condition that, only if the branch is open, i.e. $tra(t,k)=1$, $p_{top}(t,k)$ can assume values different from 0. Table 1 summarizes the parameter settings required to simulate each specific branch typology.

| | ZIL | HVDC | VSSA |
|---|---|---|---|
| $B_{min}(t,k)$, $B_{max}(t,k)$ | Set to 0 | Set to 0 | Lower and upper bound of variable susceptance |
| $p_c(t,k)$ | Set to ZIL capacity | Set to HVDC capacity | Set to 0 |
| $\Delta\delta_{max}(t,k)$ | Set to 0 | Maximum allowed phase difference | Maximum allowed phase difference |
| $\varphi_{min}(t,k)$, $\varphi_{max}(t,k)$ | Set to 0 | Set to 0 | Lower and upper bound of phase shift angle |
| $f_{max}(t,k)$ | Branch capacity | Branch capacity | Branch capacity |

**Table 1** – Parameter settings for different branch typologies.

## 3 Numerical applications

All the numerical applications presented below are implemented in MATLAB environment with CPLEX 12.9 solver, using an Intel Core i7-4770 3.4 GHz CPU with 32 GB RAM. In this section we show results obtained with the inclusion of the model in two classical DCOPF formulations, such as network constrained unit commitment (NCUC) and economic dispatch (ED). The aim is to show how the proposed model is able to account for the very diverse available remedial actions, at a low computational expense.

### 3.1 Application to NCUC

The global branch model has been applied to a seven time stamps ($N_T=7$) NCUC using the 'case6ww' 6-bus system, available in [13], as test case. The NCUC model formulated in [14, pp. 216-222] has been implemented, integrating in it the proposed branch model. The 6-bus system consists of 6 buses, 11 branches and 3 generators (connected at buses 1, 2 and 3, respectively); technical and economic data of the generating units are given in Table 2; all generating units are on-line prior to the first time stamp of the planning horizon. Regarding demands in the 7 h planning horizon, considering the load in [13] as the base load, all loads are multiplied by a coefficient, equal to 1, 1.3, 0.8, 1.1, 0.7, 0.6 and 1, respectively for each time stamp from 1 to 7. The required reserve for each time period is 10% of the aggregate load, whereas branch capacities are the same as the ones reported in [13]. The NCUC problem has been solved for five different cases, each reflecting a different flexibility of the power system: a) no remedial action is at disposal; b) branch 3-6 is a PST, with a shift angle in the (–20°, +20°) range; c) with respect to case b), branch 1-4 is a TSCS with a variable susceptance in a ±30% range around the given value in [13]; d) with respect to case c), branch 3-5 is an HVDC line, with the same capacity of the original HVAC branch; e) with respect to case d), transmission switching is allowed for each branch. Total costs and computation times of the 7 h planning horizon NCUC problem for the described cases are reported in Table 3, whereas Tables 4 and 5 show the power outputs of generating units and the power flows for all branches in case e), respectively.

| | Unit 1 | Unit 2 | Unit 3 |
|---|---|---|---|
| Minimum power output (MW) | 50 | 37.5 | 45 |
| Maximum power output (MW) | 200 | 150 | 180 |
| Ramping-down limit (MW/h) | 150 | 100 | 130 |
| Shut-down ramping limit (MW/h) | 150 | 100 | 130 |
| Ramping-up limit (MW/h) | 70 | 50 | 50 |
| Start-up ramping limit (MW/h) | 70 | 50 | 50 |
| Fixed cost ($) | 213.1 | 200 | 240 |
| Start-up cost ($) | 10 | 15 | 10 |
| Shut-down cost ($) | 0.5 | 1 | 0.3 |
| Variable cost ($/MWh) | 10 | 15 | 20 |

**Table 2** – Data of generating units.

| | Case a) | Case b) | Case c) | Case d) | Case e) |
|---|---|---|---|---|---|
| Total costs ($) | 316156.2 | 296839.7 | 296729.4 | 286082.7 | 286040 |
| Time (s) | 0.1224 | 0.1730 | 0.1314 | 0.1163 | 0.4788 |

**Table 3** – Total costs incurred in each case.

| | 1 | 2 | 3 | 4 | 5 | 6 | 7 |
|---|---|---|---|---|---|---|---|
| Unit 1 | 127 | 140 | 127 | 140 | 71 | 0 | 70 |
| Unit 2 | 83 | 133 | 41 | 91 | 76 | 126 | 140 |
| Unit 3 | 0 | 0 | 0 | 0 | 0 | 0 | 0 |

**Table 4** – Power outputs (MW) of generating units in the 7 time stamps of case e) NCUC.

| | 1-2 | 1-4 | 1-5 | 2-3 | 2-4 | 2-5 | 2-6 | 3-5 | 3-6 | 4-5 | 4-6 |
|---|---|---|---|---|---|---|---|---|---|---|---|
| 1 | 31.3 | 55.7 | 40 | 18.8 | 22.9 | 19.1 | 53.5 | 18.8 | OFF | 8.6 | 16.5 |
| 2 | 40 | 60 | 40 | 40 | 32.8 | 13.3 | 86.9 | 35.9 | 4.1 | 1.8 | OFF |
| 3 | 40 | 54.8 | 32.2 | 40 | 4.3 | 5.6 | 31.2 | 15.2 | 24.8 | 3.1 | OFF |
| 4 | 40 | 60 | 40 | 40 | 21.6 | 13.3 | 56.1 | 43.1 | -3.1 | 4.6 | 24 |
| 5 | 34.3 | OFF | 36.7 | OFF | 47.5 | 13.8 | 49 | OFF | OFF | -1.51 | OFF |
| 6 | -27.2 | 1.5 | 25.8 | 40 | 58.8 | OFF | OFF | -2 | 42 | 18.2 | OFF |
| 7 | 15.4 | 26.6 | 27.9 | 40 | 45.3 | 17.6 | 52.6 | 40 | OFF | 1.9 | 17.4 |

**Table 5** – Power flows (MW) through all 11 branches in the 7 time stamps of case e) NCUC (OFF means that the branch is switched).

Results in Table 3 show that the proposed model is able to exploit properly the power system flexibility: indeed generating costs decrease at the increasing of remedial actions available. Moreover, with the exception of transmission switching in case e), computation times are not considerably affected by the inclusion of remedial actions in the model: on the contrary, in case d) computation time is slightly less than in case a). Lastly, Table 5 shows that the model allows switching any kind of branch when needed, for instance PST in time stamps 1, 5 and 7 and TSCS in time stamp 5.

**3.2 Application to ED**

The model has been applied to single snapshots ($N_T$=1) of two network constrained ED problems, using the IEEE 118-bus system and the large scale Polish 2383-bus system as test cases. Both networks are available in [13] and in the following are referenced as base cases. Equations (9)-(11) below are added to the model (1)-(8) for the complete ED formulation:

$$\min \boldsymbol{c}^T \cdot \boldsymbol{p} \tag{9}$$

$$s.t.\ \boldsymbol{0} \le \boldsymbol{p} \le \boldsymbol{p}_{\boldsymbol{max}} \tag{10}$$

$$\boldsymbol{nex} = \boldsymbol{M}_{\boldsymbol{g}} \cdot \boldsymbol{p} - \boldsymbol{d} \tag{11}$$

Equation (9) minimizes generation costs, (10) defines the output limits of generators and (11) describes the net active power injection at each bus of the system ($\boldsymbol{M}_{\boldsymbol{g}}$ is the $N_N \times N_G$ matrix mapping generators into nodes, being $N_G$ the number of generating units in the network).

The IEEE 118-bus system consists of 118 buses, 185 branches and 19 generators; dispatch cost in the base case is 2074.4 $/h, decreasing to 1303.3 $/h in a transportation model of the network. Differently from the base case, the following remedial actions are made available: branch 81-80 is a PST, with a shift angle in the (–15°, +15°) range; branches 68-61, 26-30 and 63-64 are HVDC lines; branches 69-77, 77-82 and 89-92 are VSSAs with a variable susceptance in a ±50% range around the given value in the base case and a shift angle in the (–15°, +15°) range. Branch capacities are the same of the base case. By exploiting such remedial actions, dispatch cost becomes 1355 $/h (about 34.7% less than the base case); if TS is also allowed for all branches, the dispatch cost of the transportation model is reached by switching off branches 65-66 and 82-96. In all cases, computation times are under 450 ms.

The Polish 2383-bus system consists of 2383 buses, 2896 branches and 327 generators. In the base case the system is subdivided in six zones, of which zone 6 represents an equivalent of the foreign networks; the dispatch cost is 1,786,648 $/h, decreasing to 1,761,680 $/h with a transportation model of the network. Differently from the base case, the following remedial actions are considered: branches connecting buses of zone 6 with other buses are HVDC lines (overall 13 branches); all other branches connecting buses of different zones are VSSAs with a variable susceptance ±20% around the given value in the base case and a shift angle in the (–30°, +30°) range (overall 79 branches); branches at 400 kV voltage level, with resistance and charging susceptance equal to 0 p.u. and reactance equal to 0.0001 p.u. in the base case are considered ZILs (overall 18 branches); transmission switching is only allowed for ZILs, i.e. constraints $tra(t,k)=0$ if branch $k$ is not a ZIL are added to the formulation in (1)-(8). By exploiting such remedial actions and switching off 3 ZILs (namely branches 76-75, 99-98 and 128-127), dispatch cost becomes 1,762,127 $/h, with a 24521 $/h reduction; computation time is about 1.5 s.

**4 Conclusions**

This short communication deals with a global circuit model formulation able to represent in a DC power flow environment a wide array of branch typologies, such as variable reactance and/or variable shift angle devices, AC and HVDC lines, open branches and zero impedance lines. Reported results in case of network constrained unit commitment and economic dispatch show the potential benefit of exploiting simultaneously a wide array of remedial actions of different type. We encourage the use of the proposed model in DCOPF applications to incorporate a good level of power system flexibility at a very low computational cost.